\documentclass[12pt]{article}
\usepackage{amsmath,color,amssymb,epsfig,graphicx}
\newtheorem{theorem}{Theorem}\newtheorem{lemma}{Lemma}
\newtheorem{corollary}{Corollary}
\newtheorem{claim}{Claim}

\def\beq{\begin{equation}}\def\eeq{\end{equation}}
\def\beqn{\begin{eqnarray}}\def\eeqn{\end{eqnarray}}
\def\pont{\hspace{-6pt}{\bf.\ }}
\def\qed{\ifhmode\unskip\nobreak\fi\quad\ifmmode\Box\else$\Box$\fi}
\def\ex{{\rm ex}}

\begin{document}
\title{\bf    The linear Tur\'an number of the $k$-fan}

\pagestyle{myheadings} \markright{{\small{\sc Z.~F\"uredi and A.~Gy\'arf\'as:   On linear Tur\'an numbers}}}

\author{\bf Zolt\'an F\"uredi
    and Andr\'as Gy\'arf\'as}

\date{
Alfr\'ed R\'enyi Institute of Mathematics,
Budapest, Hungary \\
(e-mail: \texttt{z-furedi@illinois.edu, gyarfas.andras@renyi.mta.hu})
}

\footnotetext{
{Research was supported in part by grant (no. K116769) from the National Research, Development and Innovation Office – NKFIH.
The research of the first author was also supported in part by the Simons Foundation Collaboration Grant \#317487.}
 \emph{2010 Mathematics Subject Classification}:
  05D05, 05B07, 05C70.\hfill  September 14, 2017.}

\maketitle

\begin{abstract}
A hypergraph is {\em linear} if any two edges intersect in at most one vertex. For a fixed $k$-uniform family ${\cal{F}}$ of hypergraphs, the linear Tur\'an number $\ex_{\rm lin}(n,{\cal{F}})$ is the maximum number of edges in a $k$-uniform linear hypergraph $\mathcal H$ on $n$ vertices that does not contain any member of ${\cal{F}}$ as a subhypergraph.

For $k\ge 2$ the $k$-fan $F^k$ is the  $k$-uniform linear hypergraph having $k$ edges $f_1,\dots,f_k$ pairwise intersecting in the same vertex $v$ and an additional edge $g$ intersecting all $f_i$ in a vertex different from $v$. We prove the following extension of Mantel's theorem  $$\ex_{\rm lin}(n,F^k)\le {n^2 / k^2}.$$
Moreover, $|{\mathcal H}|=n^2/k^2$ holds if and only if $n\equiv 0\pmod k$ and $\mathcal H$ is a transversal design on $n$ points with $k$ groups.

We also study $\ex_{\rm lin}(n,{\cal{F}})$ where $\cal{F}$ is any subset of the three linear triple systems with four triples on at most seven points.
\end{abstract}

\section{Definitions, the linear Tur\'an number}\label{intr}

A {\em $k$-uniform hypergraph ${\cal{H}}=(V,E)$} has vertex set $V$ and its edge set $E$ consists of some $k$-element subsets  of $V$ (repeated edges are not allowed).
A hypergraph is {\em $k$-partite}
if its vertices can be partitioned into $k$ {\em groups} so that each edge has exactly one vertex from each group.

For a fixed $k$-uniform family ${\cal{F}}$, the {\em Tur\'an number $\ex(n,{\cal{F}})$} is the maximum number of edges in a $k$-uniform hypergraph $\mathcal H$ of $n$ vertices that is $\cal{F}$-free, i.e., $\mathcal H$ does not contain any member of ${\cal{F}}$ as a subhypergraph. The systematic study of Tur\'an numbers for graphs was started with the paper of Tur\'an~\cite{TU} determining $\ex(n,K_t)$.
The case $\ex(n,K_3)=\lfloor {n^2/ 4}\rfloor$ was discovered earlier~\cite{MA} and usually cited as Mantel's theorem.
Tur\'an~\cite{TU2} also initiated extensions for hypergraphs, the surveys~\cite{FUSUR}--\cite{ROPH} are covering many developments.

A hypergraph is called {\em linear} if any two edges intersect in at most one vertex. Notice that $2$-uniform linear hypergraphs are graphs. In design theory  $3$-uniform linear hypergraphs are called {\em partial triple systems} and fixed small ones are called {\em configurations}. See the seminal monograph by Colbourn and Rosa~\cite{CORO}.

For a fixed family of hypergraphs ${\cal{F}}$, the {\em linear Tur\'an number $\ex_{\rm lin}^{k}(n,{\cal{F}})$} is the maximum number of edges in a $k$-uniform {\em linear} hypergraph  $\mathcal H$ on $n$ vertices that is {\em $\cal{F}$-free}, i.e.,  $\mathcal H$ does not contain any member of ${\cal{F}}$ as a subhypergraph. The $k$-uniform {\em linear} $\cal{F}$-free hypergraphs  $\mathcal H$ on $n$ vertices with $\ex_{\rm lin}^{k}(n,{\cal{F}})$ edges are called {\em extremal hypergraphs}. The upper index $k$ is often omitted.

Linear Tur\'an numbers studied mostly implicitly. For example, the upper bound of Ruzsa and Szemer\'edi~\cite{RUSZE} on triple systems not carrying three edges on six vertices is equivalent to $\ex_{\rm lin}^3(n,T)=o(n^2)$ where $T$ is the linear triangle (a configuration $\{ 123, 345, 561\}$).
A recent remarkable work on this topic is~\cite{CGJ} where the linear Tur\'an numbers of $k$-uniform linear cycles were studied.

\subsection{$k$-fans and transversal designs}

In this paper we study the linear Tur\'an number $\ex_{\rm lin}(n,F^k)$ where $F^k$, {\it the $k$-fan},  is defined as follows.
For $k\ge 2$, $F^k$ is the linear $k$-uniform hypergraph having $k+1$ edges $f_1,\dots,f_k$ and $g$ such that $f_1,\dots,f_k$
all contain the same vertex $v$ ({\em the center}) and the additional {\em crossing edge $g$} intersects all $f_i$ in a vertex different from $v$.
Note that $F^2$ is the graph triangle and $F^3$  (sometimes called a {\em sail}) is configuration $C_{15}$ in the list of small configurations in~\cite{CORO}.

A {\em transversal design $T(n,k)$ on $n$ vertices with $k$ groups} is a $k$-partite hypergraph with groups of equal size (thus $n$ is a multiple of $k$) and each pair of vertices from different groups is covered by exactly one hyperedge.
It is well-known that transversal designs $T(n,k)$ (or equivalently $k-2$ Mutually Orthogonal Latin Squares, MOLS) exist for all $n$ when  $n> n_0(k)$ and $k$ divides $n$.
 A {\em truncated design} is obtained from a transversal design by removing one vertex (and all edges containing it).

Our main result is the following generalization of Mantel's theorem.
\begin{theorem}\pont\label{main} One has $\ex_{\rm lin}(n,F^k)\le {n^2/ k^2}$ for all $k\ge 2$. \\
The only extremal hypergraphs are the transversal designs on $n$ vertices with $k$ groups.
\end{theorem}

To determine  $\ex_{\rm lin}(n,F^k)$ exactly for all $n$ seems to be difficult.
We managed to apply the proof method of  Theorem~\ref{main} to handle the case $n\equiv -1\pmod k$.
\begin{theorem}\pont\label{var}
One has  $\ex_{\rm lin}(n,F^k)\leq m^2+m$ for $k\ge 2$, $n=k(m+1)-1$. \\\
Truncated designs  obtained from a $T(n+1,k)$  are extremal.  Further extremal hypergraphs can be found in  Theorem~\ref{3m+2}.
\end{theorem}

The traditional Tur\'an number $\ex(n,F^k)$ is determined (for fixed $k$ and $n$ large enough) by Mubayi and Pikhurko in~\cite{MP}. The extremal configuration in that case is the complete $k$-partite $k$-uniform hypergraph with almost equal parts.

\subsection{Extremal triple systems}

Note that in Theorem~\ref{var} the extremal hypergraphs are not characterized. In fact, not only the truncated designs  are extremal.
In  this subsection in Theorem~\ref{3m+2} we complete the case $k=3$, $n=3m+2$.
For $3$-uniform hypergraphs we use the term {\em triple systems} and its vertices and edges are called {\em points and triples}.

Assume that $G$ is a graph with a proper edge-coloring, i.e., the edge set of $G$ is the union of $t$ matchings (pairwise disjoint edges) $M_1,\dots,M_t$. The {\em extension} of $G$ is the linear triple system obtained by extending $V(G)$ with $t$ new points $v_1,\dots,v_t$  and defining triples by adding $v_i$ to every edge of $M_i$ for $i=1,\dots,t$.

The Wagner graph is the eight-cycle with its long diagonals. The graph $C_{5,2}$ is the blow-up of a five-cycle $v_1,\dots,v_5$ obtained by replacing its vertices with two nonadjacent vertices $p_i,q_i$ and by placing a complete bipartite graph between $\{p_i,q_i\}$ and $\{p_{i+1},q_{i+1}\}$ for $i=1,\dots,5$ $\pmod 5$.

\begin{theorem}\pont \label{3m+2} One has $\ex_{\rm lin}(n,F^3)=m^2+m$ for $n=3m+2$.  The only extremal triple systems are the following ones:

\smallskip
\noindent ${}$\enskip {\rm (\ref{3m+2}.1)}${}$\enskip
any truncated design obtained from a transversal design $T(3m+3,3)$,

\smallskip
\noindent ${}$\enskip {\rm (\ref{3m+2}.2)}${}$\enskip
  for $m=3$, the extension of a partition of the edges of the Wagner graph into three matchings of size four,

\smallskip
\noindent ${}$\enskip {\rm (\ref{3m+2}.3)}${}$\enskip
  for $m=4$, the extension of a partition of the edges of the graph $C_{5,2}$ into four matchings of size five.
\end{theorem}

\subsection{Small configurations and $(7,4)$-dense triple systems}

It leads to interesting new and classical problems if we consider $F^3$ as a member of $\mathcal A$, the family of linear triple systems with four triples on at most seven points.
Each configuration of $\mathcal A$ contains the triangle  $T$ (with triples $123, 345,156$) and an additional triple.
Apart from isomorphisms there are three possibilities to add the fourth triple.
One can add  $367$, or $246$, or $267$ and obtains either  $F^3$  (the configuration $C_{15}$ in the monograph~\cite{CORO}), or the Pasch configuration $P$ (called $C_{16}$ in~\cite{CORO}), or
  $C_{14}$, respectively, (see Figure~1).
\begin{figure}[ht]
\begin{center}
\includegraphics[width=4in]{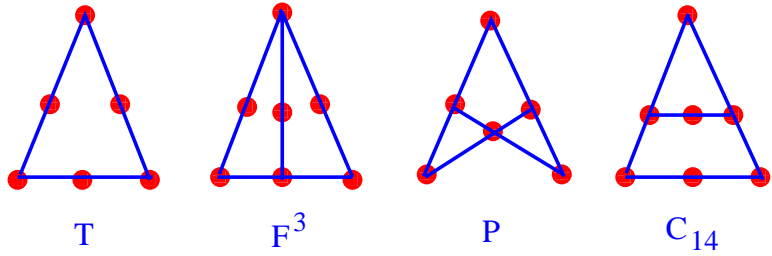}
\end{center}
\caption{$T$, $F^3$, $P$, $C_{14}$}
\label{forb}
\end{figure}

In Section~\ref{triples} we consider $\ex_{\rm lin}(n,{\cal{F}})$ for all ${\cal{F}}\subseteq {\cal{A}}=\{F^3,P,C_{14}\}$, and its relation to the notoriously difficult question (among the ones coming from~\cite{BES}): the maximum number of triples in a triple system that does not carry four triples on seven points. Let $F(7,4)$ denote the set of triple systems with four triples on at most seven points.

\begin{theorem}\pont \label{class}

\smallskip
\noindent ${}$\enskip {\rm (\ref{class}.1)}${}$\enskip
$\ex_{\rm lin}(n,P)={n(n-1)/ 6}$ \enskip for infinitely many $n$,

\smallskip
\noindent ${}$\enskip {\rm (\ref{class}.2)}${}$\enskip
$\ex_{\rm lin}(n,C_{14})={n(n-1)/ 6}$ \enskip  for infinitely many $n$,

\smallskip
\noindent ${}$\enskip {\rm (\ref{class}.3)}${}$\enskip
$\ex_{\rm lin}(n,F^3)={n^2/ 9}$\enskip  for every $n$ divisible by $3$,

\smallskip
\noindent ${}$\enskip {\rm (\ref{class}.4)}${}$\enskip
$\ex_{\rm lin}(n,\{F^3,P\})={n^2/ 9}$\enskip  for infinitely many $n$,

\smallskip
\noindent ${}$\enskip {\rm (\ref{class}.5)}${}$\enskip
$\ex_{\rm lin}(n,\{F^3,C_{14}\})={n^2/ 9}$\enskip  for infinitely many $n$,

\smallskip
\noindent ${}$\enskip {\rm (\ref{class}.6)}${}$\enskip
$\ex(n,F(7,4))-\frac{1}{2}n\le  \ex_{\rm lin}(n,\{F^3,P,C_{14}\})\le \ex_{\rm lin}(n,\{P,C_{14}\})\le  \frac{9}{2}\ex(n,F(7,4))$.
\end{theorem}

\section{$F^k$-free linear hypergraphs, Proofs} \label{kfan}

\noindent {\bf Proof of Theorem~\ref{main}. }\quad  Assume that $H$ is a $k$-uniform  $F^k$-free linear hypergraph with $n$ vertices. Let $N_v$ denote the open neighborhood of $v\in V(H)$,
$N_v:=\left(\cup_{e\ni v} e \right)\setminus \{v\}$
and set $B_v:=V(H)\setminus N_v$. Let $\Delta=\Delta(H)$ be the maximum degree of $H$ and select $v\in V(H)$ such that $d(v)=\Delta$.
Then $|N_v|=(k-1)\Delta , \,\, |B_v|=n-(k-1)\Delta. $
Since $H$ is $F^k$-free, every edge of $H$ must intersect $B_v$.
Thus $d(v)=\Delta$ implies
\beq \label{b1}|E(H)|\le \Delta|B_v|=\Delta (n-(k-1)\Delta).\eeq
On the other hand, obviously
\beq \label{b2}|E(H)|\le \frac{n\Delta}{k}. \eeq

If $\Delta \le (n/k)$ we immediately get $|E(H)|\le {n^2/ k^2}$  from (\ref{b2}).

If $\Delta > {(n/k)}$ then (\ref{b1}) and the geometric - arithmetic mean inequality imply
$$|E(H)|\le \Delta (n-(k-1)\Delta)\le \frac{1}{4}(n-(k-2)\Delta)^2<\frac{1}{4}\left(n-\frac{(k-2)n}{k}\right)^2=\frac{n^2}{k^2}.$$
Thus in both cases $|E(G)|\le {n^2/ k^2}$, as claimed.

Moreover, if equality holds then  both  estimates (\ref{b1}), (\ref{b2}) hold with equality,    thus $n=km$ and  $H$ is an $m$-regular hypergraph.
It is left to show that $|H|=m^2$ implies that $H$ is a transversal design with $k$ groups.
Since $H$ is $m$-regular,  we have $|B_v|=km-(k-1)m=m$ for every $v\in V(H)$.
\begin{claim}\pont For every $v\in V(H)$, $B_v$ is a strongly independent set, i.e.,  every edge intersects it in at most one vertex.
\end{claim}
To prove the claim, assume that $x,y\in B_v$ and there is an edge $e\in E(H)$ containing $x,y$. Then the estimate (\ref{b1}) cannot be sharp, since $e$ is counted from both $x$ and $y$. This is a contradiction, proving the Claim.

Applying the Claim for the vertices of an arbitrary edge $e=\{v_1,v_2,\dots,v_k\}$, we get the strongly independent sets $B_{v_1},\dots, B_{v_k}$. These sets must be pairwise disjoint, because if $x\in B_{v_i}\cap B_{v_j}$ then $B_{v_j}\cup \{v_i\}\subseteq B_x$, contradicting to the fact that $|B_{v_j}|=|B_x|$. Thus $V(H)$ can be partitioned into $k$ groups of size $m$, each forming a strongly independent set. The $m^2$ edges of $H$ cover $m^2{k\choose 2}$ pairs in $V(H)$ and this is equal to ${mk\choose 2}-k{m\choose 2}$, the number of pairs of $V(H)$ not covered by the groups $B_{v_1},\dots, B_{v_k}$. Thus each pair of vertices from different groups is covered exactly once, proving that $H$ is a transversal design with $k$ groups of size $m$.

Finally, it is obvious that any transversal design on $n$ vertices with $k$ groups is an $F^k$-free linear hypergraph of size $n^2/k^2$.
\qed

\bigskip

\noindent {\bf Proof of Theorem~\ref{var}. }\quad
It is obvious that any truncated design obtained from a transversal design on $n+1$ vertices with $k$ groups is an $F^k$-free linear hypergraph of size $m^2+m$.
To show that $|H|\leq m^2+m$ whenever $H$ is a $k$-uniform  $F^k$-free linear hypergraph with $n=km+k-1$ vertices we follow the argument of the proof of Theorem~\ref{main}. We use the same notations.

If $\Delta \le m$ we immediately get from (\ref{b2}) that $$|E(H)|\le \frac{((m+1)k-1)m}{k}=m^2+m-\frac{m}{k}<m^2+m.$$

If $\Delta \ge m+1$  then (\ref{b1}) and the geometric - arithmetic mean inequality imply
$$|E(H)|\le \Delta (n-(k-1)\Delta)\le \frac{1}{4}(n-(k-2)\Delta)^2\le $$ $$\le \frac{1}{4}((m+1)k-1-(k-2)(m+1))^2=m^2+m+\frac{1}{4}.$$
Thus in both cases $|E(G)|\le m^2+m$.  \qed

\bigskip

\noindent {\bf Proof of Theorem~\ref{3m+2}.}\quad  We use the notation of the proof of Theorem~\ref{var}. Assume that $H$ is an $F^3$-free linear triple system, $|V(H)|=3m+2$, $|E(H)|=m^2+m$.

 From (\ref{b2}) it follows that $\Delta \ge m+1$. The inequality (\ref{b1}) would give $|E(H)|<m^2+m$ if any vertex of $v\in V(H)$ has degree larger than $m+1$. Thus $\Delta(H)=m+1$.
Suppose $d(x)=m+1=\Delta$. Then for $B:=B_x$ one has $|B|=m$, and all vertices of $B_x$ has degree $m+1$, and $B_x$ is a strongly independent set.
It follows that $N_y=N_x$ for each $y\in B$ implying $d_H(w)=m$ for all $w\in W:=V(H)\setminus B$.

 Define the graph $G$ on vertex set $W$ and  $w_1,w_2\in W$ forming an edge in $G$ if and only if $w_1w_2v$ is a triple of $H$ for some $v\in B$. Then $G$  is $m$-regular and can be written as the union of $m$ matchings of size $m+1$. If $G$ is bipartite then must be isomorphic to $K_{m+1,m+1}$ with one matching removed, thus $H$ is a truncated design obtained from three groups of size $m+1$, the first possibility in Theorem~\ref{3m+2}.

We claim that $G$ is a triangle-free graph. Indeed, if $w_1,w_2,w_3$ are vertices of a triangle in $G$ then $H$ contains the triples
  $$e=w_1w_2v_3,\,\,  f=w_1w_3v_2,\,\,  g=w_2w_3v_1$$
for $v_i\in B_v$. Because $d_H(v_1)=m+1$, there is a triple $h=v_1w_1w_4$. Then $e,f,g,h$ form an $F^3$ with center $w_1$ and with crossing edge $g$, a contradiction.

Thus we may suppose that $G$ is a non-bipartite triangle-free graph.
A special case of a result of Andr\'asfai, Erd\H os and S\'os~\cite{AES} states that non-bipartite triangle-free graphs on $n$ vertices have minimum degree at most
$\frac{2}{5}n$. Applying this to our graph $G$, we get $m\le {2(2m+2)/ 5}$, thus $m\le 4$. There are no $m$-regular non-bipartite triangle-free graphs for $m=1,2$, but for $m=3,4$ there are: the ones in Theorem~\ref{3m+2}. Extending them with the vertices of $B$, we get the second and third possibilities in Theorem~\ref{3m+2}.  \qed

\section{Excluding small configurations of triples}\label{triples}

\subsection{Some classical transversal designs with three groups}
We refine Theorem~\ref{main} for triple systems by constructing some special transversal designs.

\begin{lemma}\pont\label{const} For every $n=6m+3$ there is a transversal design with three groups (of size $2m+1$) containing neither $F^3$ nor $P$.
\end{lemma}

\noindent {\bf Proof.} The {\em standard factorization } of the complete bipartite graph $K_{n,n}$ with partite classes $A=\{a_1,\dots,a_n\},B=\{b_1,\dots,b_n\}$ is the partition of the edge set into matchings $M_i=\{a_jb_{j+i}: 1\le j \le n\}$ for $i=0,\dots,n-1$ with indices understood $\pmod n$.  It is well-known (and it is easy to show) that in the standard factorization of the complete bipartite graph $K:=K_{2m+1,2m+1}$,  $M_i\cup M_j$ ($i\ne j$) contains no four-cycles. Considering the partite classes $A$ and $B$ of $K$ as two groups, we can add a third group $K_3=\{v_1,\dots,v_{2m+1}\}$ and for every $i=1,2,\dots,2m+1$, each edge of $M_i$ can be extended with vertex $v_i\in K_3$ to form a triple. This defines a transversal design $H$ with groups $A$, $B$, and $K_3$.
This design cannot contain $F^3$ since $F^3$ is not $3$-partite. There is no $P$ in this design either, because otherwise the symmetry of $P$ forces two points of $P$ into $K_3$ and this would force a four-cycle in the union of two $M_i$`s in $A\cup B$, a contradiction.
\qed

\begin{corollary}\pont \label{cor1} We have $\ex_{\rm lin}(n,\{F^3,P\})={n^2/ 9}$ for every $n=6m+3$ .
\end{corollary}

In a similar spirit we can also exclude from certain transversal designs the configuration $C_{14}$.

\begin{lemma}\pont\label{const2}  For every $n=3\cdot 2^m$ there is a transversal design with three groups (of size $2^m$) containing neither $F^3$ nor $C_{14}$.
\end{lemma}

\noindent {\bf Proof.} We will consider another well-known factorization of $K=[K_1,K_2]=K_{2^m,2^m}$.
The vertices of $K_1$ are labeled by $0$-$1$ sequences of length $m+1$ starting with $0$ and the vertices of $K_2$ are labeled by $0$-$1$ sequences of length $m+1$ starting with $1$.
The edges of the complete bipartite graph $[K_1,K_2]$ are labeled with the bitwise binary sum of their endpoints. Edges with the same label define the matchings of the factorization. In this factorization the union of any two matchings is partitioned into four-cycles. The transversal design $H$ is defined in the same way as in the proof of Lemma~\ref{const}, extending the $2^m$ matchings with the $2^m$ points of $K_3$. Again, there is no $F^3$ in this design. Assume that there is a $C_{14}$ in it. Apart from symmetries, there are two possible ways to distribute the vertices of $C_{14}$ in $K_1\cup K_2\cup K_3$. If two vertices of $C_{14}$ are in $K_3$, then they define an alternating path of length four in the union of two matchings in the factorization of $[K_1,K_2]$, contradiction. Otherwise three vertices of $C_{14}$ are in $K_3$, (two of degree one and one of degree two). This defines a four-cycle in $[K_1,K_2]$ whose edges are coming from (exactly) three matchings of the factorization, a contradiction again. \qed

\begin{corollary}\pont\label{cor2} We have $\ex_{\rm lin}(n,\{F^3,C_{14}\})={n^2/ 9}$  for every $n=3\cdot 2^m$.
\end{corollary}

\subsection{Proof of Theorem~\ref{class}}

$\ex_{\rm lin}(n,\emptyset)={n(n-1)/ 6}$  for infinitely many $n$ follows trivially: every Steiner triple system is an extremal hypergraph.
For $\ex_{\rm lin}(n,P)$ and for $\ex_{\rm lin}(n,C_{14})$ the extremal hypergraphs are the Steiner triple systems avoiding $P$ and the ones avoiding $C_{14}$.
This is a well studied area in the theory of Steiner triple systems, see~\cite{Br,CORO}.

The statement $\ex_{\rm lin}(n,F^3)={n^2/ 9}$ (whenever $n$ is divisible by 3) is a special case of our main result, Theorem~\ref{main}, and
$\ex_{\rm lin}(n,\{F^3,P\})={n^2/ 9}$  and $\ex_{\rm lin}(n,\{F^3,C_{14}\})={n^2/ 9}$  for infinitely many $n$ follows from Corollaries~\ref{cor1} and~\ref{cor2}, respectively.

Let $W$ denote the triple system with  two triples intersecting in two vertices.
The inequalities
\begin{align*}
\ex(n,F(7,4))-n/2&\le \ex_{\rm lin}(n,F(7,4))
\\ &= \ex_{\rm lin}(n,\{F^3,P,C_{14}\})\le \ex_{\rm lin}(n,\{P,C_{14}\})\\
\le \frac{9}{2} \ex_{\rm lin}(n,\{F^3,P,C_{14}\})&= \frac{9}{2} \ex(n,\{F^3,P,C_{14},W\})= \frac{9}{2}\ex(n,F(7,4))
    \end{align*}
showing that $\ex_{\rm lin}(n,\{P,C_{14}\})$ and $\ex_{\rm lin}(n,\{F^3,P,C_{14}\})$ have the same order of magnitude as $\ex(n,F(7,4))$, follow from the following argument.

Assume that $H$ is a triple system not containing four triples within seven points. The first inequality follows from the fact that removing the triples from the maximum number of pairwise vertex disjoint copies of $W$'s, at most $n/2$ triples are removed and the remaining triple system is linear. The next equality holds since the linear configurations carrying four triples on seven points are $F^3$, $P$, and $C_{14}$.
The next inequality is obvious from monotonicity of $\ex_{\rm lin}(n,{\cal{F}})$.
The inequality $\ex_{\rm lin}(n,\{P,C_{14}\})\le \frac{9}{2} \ex_{\rm lin}(n,\{F^3,P,C_{14}\})$ follows by the standard argument that from every triple system one can select a $3$-partite subsystem containing at least $\frac{2}{9}$ fraction of the triples. Indeed, applying this to a triple system containing no $P$ nor $C_{14}$, we get a subsystem that does not contain $F^3$ either (since $F^3$ is not $3$-partite). The next equality holds because linear triple systems are the ones that do not contain $W$. Finally the last equality follows from checking that every member of $F(7,4)$ contains at least one of $\{F^3,P,C_{14},W\}$. \qed

\eject

\end{document}